\documentclass[12pt]{article}
\usepackage{amssymb,pb-diagram,CD,papersaver6}

\newcommand{\PSbox}[3]{\mbox{\rule{0in}{#3}\includegraphics{#1}\hspace{#2}}} 
\newtheorem{Theorem}{Theorem}[section]
\newtheorem{Definition}{Definition}[section]

\newtheorem{Corollary}[Theorem]{Corollary}

\newtheorem{Fact}[Theorem]{Fact}

\newtheorem{Claim}[Theorem]{Claim}

\newcommand{\lra}{\longrightarrow}

\def\sqr#1#2{{\vcenter{\vbox{\hrule  height.#2pt
        \hbox{\vrule width.#2pt height#1pt \kern#1pt \vrule width.#2pt}
        \hrule height.#2pt}}}}
\def\Box{\sqr66}   
\let\epsilon=\varepsilon
\def\){ \right) }
\def\({ \left( }
\def\[{ \left[ }
\def\]{ \right] }
\def\<{ \langle }
\def\>{ \rangle }
\let\ljunk=\{
\let\rjunk=\}
\def\{{\left\ljunk}
\def\}{\right\rjunk}
\def\p{\partial}

\def\conf{{\mathrm c}{\mathrm o}{\mathrm n}{\mathrm f}}
\def\Riem{{\cal R}{\mathrm i}{\mathrm e}{\mathrm m}}
\def\dist{{\mathrm d}{\mathrm i}{\mathrm s}{\mathrm t}}

\def\cyl{{\mathrm c}{\mathrm y}{\mathrm l}}
\def\can{{\mathrm c}{\mathrm a}{\mathrm n}}

\newcommand{\diff}{{\mathrm d}{\mathrm i}{\mathrm f}{\mathrm f}}
\newcommand{\Diff}{{\mathrm D}{\mathrm i}{\mathrm f}{\mathrm f}}
\newcommand{\Spin}{{\mathrm S}{\mathrm p}{\mathrm i}{\mathrm n}}
\newcommand{\SO}{{\mathrm S}{\mathrm O}}
\newcommand{\PPos}{{\mathfrak P}{\mathfrak o}{\mathfrak s}}

\def\Vol{{\mathrm V}{\mathrm o}{\mathrm l}}

\newcommand{\Pos}{{\cal P}{\mathrm o}{\mathrm s}}

\newcommand{\R}{{\mathbf R}}

\newcommand{\Z}{{\mathbf Z}}

\begin{document} 

\title{Manifolds of Positive Scalar Curvature and\\
Conformal Cobordism Theory}
\author{Kazuo Akutagawa\thanks{\ \ Partialy supported by the 
Grants-in-Aid for Scientific Research, The Ministry of Education,
Science, Sports and Culture, Japan, No. 09640102
\protect\\
email: smkacta\atsign ipc.shizuoka.ac.jp
}
\and
Boris Botvinnik\thanks{
\ \ Partially supported by SFB 478, M\"unster
\protect\\
e-mail: botvinn\atsign math.uoregon.edu
}}
\maketitle
\vspace{-1.5em}
\begin{abstract}
We study here compact manifolds with positive scalar curvature
metrics.  We use the relative Yamabe invariant from \cite{AB} to define
the conformal cobordism relation on the category of such manifolds. We
prove that corresponding conformal cobordism groups
$\Pos_n^{\conf}(\gamma)$ are isomorphic to the cobordism groups
$\Pos_n(\gamma)$ defined topologically by S. Stolz in
\cite{Stolz1}. As a corollary we show that the conformal concordance
of positive scalar curvature metrics coincides with the standard
concordance relation. Our main technical tools came from the analysis
and conformal geometry.
\end{abstract}
\vspace{-2.5em}
\section{Introduction}\label{introduction}
{\bf \ref{introduction}.1. Motivation.} There are two competing
approaches in the study of manifolds admitting a metric of positive
scalar curvature. 
\vspace{2mm}

\noindent
The first approach is developed within conformal geometry and
analysis, and the second one unconventionally may be called
``topological'' (where the $Spin$-geometry and the Dirac operator
methods are combined with the differential topology and some homotopy
theory). There are recent detailed surveys presenting a current state
of affairs in the subject, given by M. Gromov \cite{G1}, and
J. Rosenberg \& S. Stolz \cite{RS}. It is emphasized in \cite{G1},
that the conformal geometry technique (which, perhaps, includes the
minimal surface method) has certain advantages over the topological
methods since it does not require $Spin$ structure and, in some
respect, even completeness of a manifold. ``On the other hand,
whenever the Dirac method applies it delivers finer geometrical (and
topological) information although in no serious case the results of
one method may be completely recaptured by the other.''$^*$\footnote{\ \ $\!^*$We quote \cite[p. 45]{G1}}
\vspace{2mm}

\noindent
The goal of this paper is to establish one particular link between the
topological and conformal approaches, where the resulting object is,
indeed, the same. Namely, we show that the cobordism groups of
manifolds with positive scalar curvature metrics, delivered by
topological means and by means of conformal geometry, coincide.
\vspace{2mm}

\noindent
{\bf \ref{introduction}.2. Restictions.}  We restrict here our
attention to the oriented smooth manifolds. There are also the
dimensional restrictions: all topological constructions work well
starting with dimension five for closed manifolds (and six for
manifolds with boundary). The conformal geometry gives the dimensional
restiction at least two (for closed manifolds) and at least three
otherwise. We use abbreviation ``psc'' for positive scalar curvature.
\vspace{2mm}

\noindent
{\bf \ref{introduction}.3. Topological psc-cobordism.}  Let
$(M_0,g_0)$, $(M_1,g_1)$ be compact manifolds with psc-metrics $g_0$
and $g_1$.  Then $(M_0,g_0)$, $(M_1,g_1)$ are {\sl psc-cobordant} if
there exists a Riemannian manifold $(W,\bar{g})$, $\p W = M_0\sqcup
(-M_1)$, so that
$$
R_{\bar{g}} > 0, \ \ \ 
\bar{g}|_{M_j}=g_j, \ \ \ \ \mbox{and $\bar{g}=g_j +dt^2$ near the
boundary $\p W= M_0\sqcup M_1$ for $j=0,1$.}
$$
We emphasize the importance of the condition that the metric $\bar{g}$
must be a product metric near the boundary. In the case of $Spin$
manifolds it gives, in particular, that the Dirac operator with the
Atiyah-Patodi-Singer boundary conditions is well-defined. In fact, as
it was showed by S. Stolz \cite{Stolz1}, each given manifold $M$
(admitting a psc-metric, and not necessarily $Spin$) lives in a
specific cobordism category $\PPos(\gamma)$, $\gamma=(\pi, w,
\hat{\pi})$, determined by the fundametal group $\pi=\pi_1(M)$ and the
Stiefel-Whitney classes $w_1(M)$, $w_2(M)$. In particular, $w : \pi
\rightarrow \Z_2$ is nothing but the orientation character given by
$w_1(M)$.  We say that the stucture $\gamma=(\pi, w, \hat{\pi})$ is
{\sl oriented} if $w =0$. Let $\Pos_n(\gamma)$ be the corresponding
psc-cobordism groups.
\vspace{2mm}

\noindent
{\bf \ref{introduction}.4. Conformal cobordism.} Let $(M,g)$ be a
compact manifold with psc-metric $g$. In the conformal geometry world,
it means that the conformal class $C=[g]$ is such that the Yamabe
constant $Y_C(M)>0$.  We call such a conformal class $C$ {\sl positive}.
Let ${\cal C}^+(M)$ be the space of positive conformal classes. We
call a pair $(M, C)$ with $C\in {\cal C}^+(M)$ a {\sl positive
conformal manifold}.  Now let $W$ be a compact smooth manifold with
boundary, $\p W = M\neq \emptyset$, and let $C$ be a conformal class
on $M$. Let $\bar{C}$ be a conformal class on $W$.  We
write $\p \bar{C} = C$ if the conformal class $\bar{C}$ is such that
$\bar{C}|_{M} = C$. We defined in \cite{AB} the {\sl relative Yamabe
constant} $Y_{\bar{C}}(W, M; C)$ and the {\sl relative Yamabe
invariant}
$$
Y(W,M; C)= \sup_{\bar{C}, \p \bar{C} = C} Y_{\bar{C}}(W,M;C).
$$
We emphasize that in order to define 
$Y_{\bar{C}}(W, M; C)$, we use the subclass $\bar{C}^0\subset \bar{C}$
of metrics with {\sl zero mean curvature} along $M$ (see
\cite{AB} and Section \ref{relative}).
\vspace{2mm}

\noindent
Positive conformal manifolds $(M_0,C_0)$, $(M_1,C_1)$ are {\sl
conformally cobordant} if there exists a smooth cobordism $W$ with 
boundary $\p W = M_0\sqcup(-M_1)$, and such that the relative Yamabe
invariant $ Y (W,M_0\sqcup M_1; C_0\sqcup C_1)>0.  $ We proved in
\cite{AB} that the conformal cobordism is an equivalence
relation.  We also incorporate the above oriented $\gamma$-structure
into this cobordism equivalence to define the conformal cobordism
groups $\Pos^{\conf}_n(\gamma)$ of positive conformal manifolds
equipped with a $\gamma$-structure. Clearly there is a natural
homomorphism $\Pos_n(\gamma)\lra \Pos^{\conf}_n(\gamma)$ given by
taking conformal classes of corresponding metrics.
\vspace{2mm}

\noindent
{\bf Remark.} Perhaps, it is important to emphasize the major
difference between the above cobordism relations.  Firstly, it is in
the boundary conditions: product metric near the boundary versus
vanishing of the mean curvature along the boundary. Secondly, let
$\bar{C}$ be a conformal class on $W$ such that the relative Yamabe
constant $ Y_{\bar{C}} (W,M_0\sqcup M_1; C_0\sqcup C_1)$ is positive:
such a conformal class $\bar{C}$ exists if the relative Yamabe
invariant $Y (W,M_0\sqcup M_1; C_0\sqcup C_1)>0$. Then a metric
$\bar{g}\in \bar{C}$ (which restricts to given psc-metrics $g_0$ and
$g_1$ on the boundary and even is a product metric near the boundary)
may not have, in general, positive scalar curvature.
\vspace{2mm}

\noindent
{\bf \ref{introduction}.5. Main results.} 
\vspace{2mm}

\noindent
{\bf Theorem A.}  {\sl Let $\gamma$ be an oriented structure, and
$n\geq 5$. Then the conformal cobordism groups
$\Pos^{\conf}_n(\gamma)$ are naturally isomorphic to the psc-cobordism
groups $\Pos_n(\gamma)$.}
\vspace{3mm}

\noindent
Recall that in the conformal world the classic Yamabe invariant $Y(M)$
gives very simple answer on the existence of psc-metric.  Indeed:
\vspace{2mm}

\noindent
$\bullet$ {\sl Let $M$ be a closed oriented manifold with $\dim M \geq
2$. Then the Yamabe invariant $Y(M)>0 $ if and only if there exists a
psc-metric on $M$.}
\vspace{2mm}

\noindent
The relative Yamabe invariant has a similar property (where the
manifolds below are oriented).
\vspace{3mm}

\noindent
{\bf Corollary B.} {\sl Let $M = \p W$, $\dim M \geq 2$, and $g$ be a
psc-metric on $M$.  Then the relative Yamabe invariant $Y(W,M;[g])>0
$ if and only if the metric $g$ may be extended to a psc-metric
$\bar{g}$ on $W$, so that $\bar{g}$ is a product metric near the
boundary.}
\vspace{3mm}

\noindent
On the topological side of this story, S. Stolz also defines the
relative cobordism groups $R_n(\gamma)$, \cite{Stolz1} (see
\cite{Hajduk} for the simply connected case).  S. Stolz proves that
the cobordism groups $R_n(\gamma)$ are the {\sl actual obstruction
groups} for the existence of psc-metrics (see \cite[Theorem
1.1]{Stolz1}). The groups $\Pos_n(\gamma)$, $R_n(\gamma)$ and the
regular cobordism groups $\Omega_n(\gamma)$ (of manifolds carrying
$\gamma$-structure) fit together into the exact sequence
\begin{equation}\label{e-1}
\cdots\rightarrow R_{n+1} (\gamma)\lra \Pos_n(\gamma) \lra
\Omega_n(\gamma) \lra R_n(\gamma)\lra \Pos_{n-1}(\gamma) \rightarrow
\cdots
\end{equation}
In the case of simply connected $Spin$ manifolds, $\Omega_n(\gamma)=
\Omega^{\Spin}_n$. We define the conformal ``relatives'' to
$R_n(\gamma)$ (the cobordism groups $R_n^{\conf}(\gamma)$) for
oriented $\gamma$-structures, so that there is the exact sequence
$$
\cdots\rightarrow R^{\conf}_{n+1} (\gamma)\lra \Pos^{\conf}_n(\gamma)
\lra \Omega_n(\gamma) \lra R^{\conf}_n(\gamma)\lra
\Pos_{n-1}^{\conf}(\gamma) \rightarrow \cdots
$$
which turns out to be isomorphic to (\ref{e-1}). In particular, we have
\vspace{3mm}

\noindent
{\bf Corollary C.} {\sl Let $\gamma$ be an oriented structure, and
$n\geq 6$. Then the conformal cobordism groups $R^{\conf}_n(\gamma)$
are naturally isomorphic to the psc-cobordism groups $R_n(\gamma)$.}
\vspace{2mm}

\noindent
{\bf \ref{introduction}.6. Concordance and conformal concordance of
psc-metrics.} Recall that two psc-metrics $g_0$, $g_1$ on $M$ are {\sl
psc-concordant} if there exists a psc-metric $\bar{g}$ on a cylinder
$M\times [\ell_0,\ell_1]$ (for some $\ell_0< \ell_1$) so that
$$
\bar{g}|_{M\times \{\ell_j\}}=g_j, \ \ \ \ \mbox{and $\bar{g}=g_j +dt^2$
near the boundary $M\times \{\ell_j\}$ for $j=0,1$.}
$$
Two positive conformal classes $C_0,C_1\in {\cal C}^+(M)$ are {\sl
conformally concordant} if the Yamabe invariant 
$$
Y(M\times [0, 1], M\times \{0,1\}; C_0\sqcup C_1)>0.
$$ 
We proved in \cite{AB} that conformal concordance is an equivalence
relation. Clearly the psc-concordance implies the conformal
concordance. 
\vspace{3mm}

\noindent
{\bf Corollary D.} {\sl Let $M$ be an oriented manifold with $\dim
M\geq 2$, and let $g_0, g_1$ be psc-metrics on $M$ such that the
conformal classes $C_0=[g_0]$ and $C_1=[g_1]$ are conformally
concordant. Then the metrics $g_0$, $g_1$ are psc-concordant.}
\vspace{3mm}

\noindent
{\bf \ref{introduction}.7. Important remark.} Unfortunately the
results of this paper do not allow to compute the conformal cobordism
groups $\Pos_n^{\conf}(\gamma)$ and $R_n^{\conf}(\gamma)$. However, in
our view, these results open up some new possiblities which hopefully
will be explored by open-minded geometers and topologists. We discuss
this in Section \ref{disscussion}.
\vspace{2mm}

\noindent
{\bf \ref{introduction}.8. Organization of the paper.} We review
necessary constructions and facts on the conformal geometry in Section
\ref{relative}. We state our main technical result in Section
\ref{outline} and outline key points of its proof.  We give this proof
in Section \ref{cobordism}.  We review some topological constructions
and finish the proofs in Section \ref{topology}. Finally we discuss
some open problems in Section \ref{disscussion}.
\vspace{2mm}

\noindent
{\bf \ref{introduction}.9. Acknowledgments.}  Both authors would like
to acknowledge partial financial support provided by the Department of
Mathematics at the University of Oregon, SFB 478 -- Geometrische
Strukturen in der Mathematik, 
and the Grants-in-Aid for
Scientific Research, Japan. We would like also to thank Michael
Joachim, Wolfgang L\"uck, and Thomas Schick for useful discussions and
warm hospitality during our visit at M\"unster.

\section{Some conformal geometry}\label{relative}
{\bf \ref{relative}.1. General setting.} Let $W$ be a compact smooth
manifold with boundary, $\p W = M\neq \emptyset$, and $n=\dim W \geq
3$. We always assume that all manifolds are oriented, and the
orientation on $W$ is compatible with the orientaion on its boundary
$\p W$.
\vspace{2mm}

\noindent
Let $C$ be a conformal class of metrics on $M$, and $\Riem(W)$ is the
space of all Riemannian metrics on $W$.  For a metric $\bar{g}\in
\Riem(W)$ we denote $H_{\bar{g}}$ the mean curvature along the
boundary $\p W = M$. We denote ${\cal C}(M)$ and ${\cal C}(W)$ the
space of conformal classes on $M$ and $W$ respectively.  Let $C\in
{\cal C}(M)$ $\bar{C}\in {\cal C}(W)$.  We say that $C$ is the {\sl
boundary of $\bar{C}$} or $\bar{C}$ is a {\it coboundary of $C$} if
$\bar{C}|_{M} = C$.  We use notation $\p \bar{C}=C$ in this case. Then
a pair of conformal classes $(\bar{C},C)$ is a {\it conformal class on
$(W,M)$ if $\p \bar{C}=C$.}  We denote ${\cal C}(W,M)$ the space of
pairs of conformal classes.  Let $(\bar{C},C)\in {\cal C}(W,M)$.  For
each pair of conformal classes $(\bar{C},C)\in {\cal C}(W,M)$ we
consider the conformal subclass $\bar{C}^0 \subset \bar{C}$ defined as
$$
\bar{C}^0 = \{\bar{g}\in \bar{C} \ | \ H_{\bar{g}}=0 \ \}.
$$
We call $\bar{C}^0 \subset \bar{C}$ the {\sl normalized conformal
class}. Let ${\cal C}^0(W,M)$ be the space of pairs $(\bar{C}^0,C)$,
so that $\bar{C}^0\subset \bar{C}$ as above, and $(\bar{C},C)\in {\cal
C}(W,M)$. In fact, it is easy to see that for any conformal class
$\bar{C}\in {\cal C}(W)$ the subclass $\bar{C}^0$ is not empty (see
\cite[formula (1.4)]{Escobar1}). Thus there is a natural bijection
between the spaces ${\cal C}^0(W,M)$ and ${\cal C}(W,M)$.  Let
$\bar{g}\in \bar{C}^0$ be a metric. Then $\bar{C}^0$ may be described
as follows:
$$
\bar{C}^0 = \{ u^{{4\over n-2}} \bar{g} \ | \ u\in C^{\infty}_+(W) \
\mbox{such that} \ \ \p_{\nu}u = 0 \ \ \mbox{along} \ M \ \}.
$$
Here $\nu$ is a normal unit (inward) vector field along the boundary, and 
$C^{\infty}_+(W)$ is the space of positive smooth functions on $W$.
\vspace{2mm}

\noindent
{\bf \ref{relative}.2. The Einstein-Hilbert functional.} Let $C\in
{\cal C}(M)$ be given. We define the following subspaces of metrics:
$$
\begin{array}{lcl}
\Riem_C(W,M) & = & \{ \bar{g}\in \Riem(W) \ | \  \p [\bar{g}]=C \},
\\
\\
\Riem_C^0(W,M) & = & \{ \bar{g}\in \Riem_C(W) \ | \   H_{\bar{g}} = 0 \}.
\end{array}
$$
The normalized Einstein-Hilbert functional $I : \Riem_C^0(W,M)
\longrightarrow \R$ is given by
$$
I(\bar{g})= {\int_W R_{\bar{g}} d\sigma_{\bar{g}}\over
\Vol_{\bar{g}}(W)^{{n-2\over n}}},
$$
where $R_{\bar{g}}$ is the scalar curvature, and $d\sigma_{\bar{g}}$
is the volume element. The following fact is analogous to the classic
theorem on the  Einstein-Hilbert functional.
\begin{Theorem}\label{Th1}{\rm (\cite[Theorem 1.1]{AB})}
Critical points of the functional $I$ on the space of metrics
$\Riem_C^0(W,M)$ coincide with the set of Einstein metrics $\bar{g}$
on $W$ with $\p [\bar{g}]=C$, and $H_{\bar{g}} = 0$.
\end{Theorem}

\noindent
{\bf \ref{relative}.3. Relative Yamabe invariant.} Let $(\bar{C},C)\in
{\cal C}(W,M)$.  The {\it relative Yamabe costant of} $(\bar{C},C)$ is
defined as
$$
Y_{\bar{C}}(W,M;C) = \inf_{\bar{g}\in \bar{C}^0} I(\bar{g}).
$$
{\bf Remark.} The relative Yamabe constant $Y_{\bar{C}}(W,M;C)$ is
related to the Yamabe problem on a manifold with boundary, which was
solved by P. Cherrier \cite{Cherrier} and J. Escobar \cite{Escobar1}
under some restrictions. In fact, in a generic case there is a
relative {\sl Yamabe metric} $\check{g}\in \bar{C}$ with
$H_{\check{g}}=0$ and constant scalar curvature $R_{\check{g}}=
Y_{\bar{C}}(W,M;C)\cdot \Vol_{\check{g}}(W)^{-{2\over n}}$ (see
\cite{Escobar1}, \cite{AB} for more details).
\vspace{2mm}

\noindent
The {\sl relative Yamabe invariant} with respect to a conformal class
$C\in {\cal C}(M)$ is defined as:
$$
Y(W,M;C) = \sup_{\bar{C}, \p \bar{C} = C}
Y_{\bar{C}}(W,M;C).
$$
The relative Yamabe invariant $Y(W,M;C)$ has several important
properties analogous to the corresponding properties of the classic
Yamabe invariant (see \cite{AB} for details).
\vspace{2mm}

\noindent
{\bf \ref{relative}.4. Approximation theorem.} One notices that the
minimal boundary condition is crucial to define the relative Yamabe
constant. In general, it is rather delicate problem to approximate
given conformal class $\bar{C}$ on a manifold with boundary by such
conformal classes which contain a product metric near the boundary
(see \cite{Kobayashi1}, \cite{AB}). The minimal boundary condition is
crucial to prove the following approximation result.
\begin{Theorem}\label{main-trick}{\rm 
(\cite[Theorem 4.6]{AB})}
Let $W$ be a manifold with boundary $\p W= M$, $C\in {\cal C}(M)$.
Let $\bar{g}\in \Riem_C^0(W,M)$ be a metric. Let $g= \bar{g}|_M$, and
$A_{\bar{g}}$ be the second fundamental form of $M=\p W$. There exists
a family of metrics $\tilde{g}_{\delta}\in \Riem_C^0(W,M)$ such that
\begin{enumerate}
\item[{\bf (i)}] $\tilde{g}_{\delta} \rightarrow \bar{g}$ in
the $C^0$-topology on $W$ {\rm (}as $\delta  \rightarrow 0${\rm )},
\item[{\bf (ii)}] $R_{\tilde{g}_{\delta}} \rightarrow R_{\bar{g}}$ in
the $C^0$-topology on $W$ {\rm (}as $\delta  \rightarrow 0${\rm )},
\item[{\bf (iii)}] $\tilde{g}_{\delta}$ conformally equivalent to
the metric $g + dr^2$ on 				
$U_{\epsilon(\delta)}(M,\bar{g})$,
\item[{\bf (iv)}]  $\tilde{g}_{\delta} \equiv \bar{g}$ on $W\setminus 
U_{\delta}(M,\bar{g})$.
\end{enumerate}
\end{Theorem}
In terms of the relative Yamabe constant, 
Theorem \ref{main-trick} gives the following conclusion:
\begin{Corollary}\label{Cor-trick}{\rm 
(\cite[Theorem 2.1]{AB})}
For any $\bar{C}\in {\cal C}_C(W,M)$, and any $\epsilon>0$ there
exist a conformal class $\tilde{C}\in {\cal C}_C(W,M)$, and a
metric $\tilde{g}\in \tilde{C}^{0}$, such that
\begin{equation}
\{
\begin{array}{l}
\bar{C} \ \mbox{and} \ \tilde{C} \ \mbox{are $C^0$-close conformal classes}
\\
|Y_{\bar{C}}(W,M;C) - Y_{\tilde{C}}(W,M;C)| < \epsilon
\\
\tilde{g} \sim g + dr^2 \ \ \mbox{(conformally equivalent near $M$),}
\end{array}
\right.
\end{equation}
where $C=\p\tilde{C}$ and $g = \bar{g}|_M$. 
\end{Corollary}

\noindent
{\bf \ref{relative}.5. Conformal cobordism.} We call a conformal class
$C\in {\cal C}(M)$ {\sl positive} if the Yamabe constant $Y_C(M)>0$.
Let ${\cal C}^+(M)\subset {\cal C}(M)$ be the space of all positive
conformal classes. A pair $(M,C)$ with $C\in {\cal C}^+(M)$ is called
a {\sl positive conformal manifold}. Recall that positive conformal
manifolds $(M_0,C_0)$, $(M_1,C_1)$ are {\sl conformally cobordant} if
there exists a smooth cobordism $W$ with the boundary $\p W =
M_0\sqcup(-M_1)$, and such that the relative Yamabe invariant $ Y
(W,M_0\sqcup M_1; C_0\sqcup C_1)>0$. We proved in \cite{AB} that the
conformal cobordism is an equivalence relation.
\vspace{2mm}

\noindent
{\bf \ref{relative}.5. Cylindrical manifolds.} It is convenient for us
to use a general concept of cylindrical manifolds. Let $Z$ be a
compact, closed smooth manifold, $\dim Z= n-1$. In general, $Z$ may
have several connective components; we let
$$
Z = \bigsqcup_{j=1}^m Z_j, \ \ \ \mbox{where each $Z_j$ is connected.}
$$
Let $\Riem(Z)$ be the space of Riemannian metrics on the manifold
$Z$. We let $h\in \Riem(Z)$ to be fixed.
\begin{Definition}\label{Def1}
{\rm Let $(X, \bar{g})$ be a complete Riemannian manifold, $\dim X=n$. We
call $(X, \bar{g})$ a {\it cylindrical manifold modeled by $(Z,h)$} if
there exists a compact smooth manifold $W$ with non-empty boundary $\p
W= Z\sqcup M$ such that
$$
\{
\begin{array}{l}
X \stackrel{\diff}{\cong} W\cup_Z \( Z \times [0,\infty)\) \ \ \ 
\mbox{where $\p W \supset Z$ is identified with $Z \times\{0\}\subset 
Z \times [0,\infty)$,}
\\
\\
\bar{g}(z,t) = h(z) +dt^2  \ \ \ 
\mbox{on $Z\times [1,\infty)$ with $(z,t)$ coordinates on 
$Z\times [1,\infty)$}
\end{array}
\right.
$$ 
(see Fig. \ref{relative}.1). The metric $\bar{g}$ is called a {\sl
cylindrical metric on $X$}.}
\end{Definition}
We define the space of cylindrical metrics on $X$:
$$
\Riem^{\cyl}(X) = \{ \bar{g}\in \Riem(X) \ | \ \bar{g} \ \mbox{is
cylidrical as in Definition \ref{Def1}}\}
$$

\hspace*{20mm}
\PSbox{cob1.pstex}{10mm}{44mm}  
\begin{picture}(0,1)
\put(70,110){{\small $Z\!\times\!\!\{0\}$}}
\put(110,110){{\small $Z\!\times\!\!\{1\}$}}
\put(170,110){{\small $Z\!\times\!\! [1,\infty)$}}
\put(83,-9){{\small $0$}}
\put(120,-9){{\small $1$}}
\put(295,0){{\small $\R$}}
\put(-30,103){{\small $M$}}
\put(17,78){{\small $W$}}
\end{picture}
\vspace{4mm}

\centerline{{\small {\bf Fig. \ref{relative}.1.} A cylindrical manifold $X$}}
\vspace{2mm}

\noindent
We define the space of {\sl cylindrical conformal classes} on $X$ as 
$$
{\cal C}^{\cyl}(X) = \{ [\bar{g}] \ | \ \bar{g} \in \Riem^{\cyl}(X)
\}.  
$$
{\bf Remark.} The category of cylindrical manifolds is well-suited for
the conformal geometry. In particular, there are well-defined {\sl
cylindrical} Yamabe constant and Yamabe invariant. The authors plan to
explore this in different paper.

\section{Main Theorem: outline of the proof}\label{outline}
{\bf \ref{outline}.1. Setting.} Let $W$ be a compact smooth manifold
with $\dim W =n\ge 3 $, and $\p W= Z\sqcup M \neq \emptyset$. Let
$C\in {\cal C}^{+}(Z)$ be a positive conformal class, and $C^{\prime}
\in {\cal C}(M)$ a conformal class. (Here $C^{\prime}$ may not be
positive, in general.)  We choose a metric $h\in C\sqcup C^{\prime}$
with $R_h>0$. We assume that it is given a conformal class
$\tilde{C}\in {\cal C}(W)$ with $\p \tilde{C}= C\sqcup C^{\prime}$,
and that $Y_{\tilde{C}}(W,Z\sqcup M;C\sqcup C^{\prime})>0$. We use
Theorem \ref{main-trick} and \cite[Theorem 5.1]{AB} to choose a
conformal class
$$
\bar{C}\in  {\cal C}(W\cup_Z (Z\times [0,1]))= {\cal C}(X(1)) 
$$
on the manifold $X(1)= W\cup_Z(Z\times [0,1])$ satisfying the
following properties:
\begin{enumerate}
\item[{\bf (1)}] the restriction $\bar{C}|_W$ is a small pertubation
of $\tilde{C}$ near $Z\times\{0\}$;
\item[{\bf (2)}] there is a metric $\bar{g}\in \bar{C}$ such that
$\bar{g}= h+dt^2$ near $Z\times \{1\}$;
\item[{\bf (3)}] $Y_{\bar{C}}(X(1),Z\sqcup M;C\sqcup C^{\prime})>0$.
\end{enumerate}
We extend ``cylindrically'' the metric $\bar{g}$ and the conformal
class $\bar{C}\in {\cal C}(X(1))$ to the cylindrical manifold
$$
X = W\cup_Z\(Z\times [0,\infty)\) = X(1)\cup_Z \(Z\times [1,\infty)\).
$$
We denote those extensions also $\bar{g}\in \Riem^{\cyl}(X)$, and
$\bar{C}\in {\cal C}^{\cyl}(X)$. The resulting manifold $X$ is a
cylindrical manifold modeled by $(Z,h)$. 
Thus we have 
$$
\{
\begin{array}{l}
Y_{\bar{C}|_{X(1)}}(X(1), Z\sqcup M;C\sqcup C^{\prime}) > 0 ,
\\
\bar{g} = h +dt^2 \ \ 
\mbox{on $Z\times (1-\epsilon,\infty)$ for some $\epsilon>0$.}
\end{array}
\right.
$$
We define $X(\ell)= W\cup_Z\( Z\times [0,\ell]\)$ for
$\ell\geq 1$.
\vspace{2mm}

\noindent
The following theorem is the main technical result in this section. We
use the above notations in this theorem.

\begin{Theorem}\label{main-new}
Let $W$ be a compact smooth manifold, $\dim W \geq 3$, with $\p W = Z
\sqcup M$, $Z\neq \emptyset$. Let $C\in {\cal C}^+(Z)$, $C^{\prime}\in
{\cal C}(M)$, and let $h\in C$ be a given metric with $R_h>0$.
\vspace{2mm}

\noindent
Let $\tilde{C}\in {\cal C}(W)$ be a conformal class with $\p
\tilde{C}= C\sqcup C^{\prime}$, such that the relative Yamabe constant
$Y_{\tilde{C}}(W,Z\sqcup M ;C\sqcup C^{\prime})>0$. Let $X$ be the
above cylindrical manifold modeled by $(Z,h)$.
\vspace{2mm}

\noindent
Then there exist a constant $L>>1$, a conformal class $\bar{C}\in
{\cal C}^{\cyl}(X(L))$ with $\p \bar{C}=C\sqcup C^{\prime}$, and a
metric $\hat{g}\in \bar{C}^0$, such that
\begin{equation}\label{l12}
\{
\begin{array}{ll}
R_{\hat{g}} > 0 & \mbox{on $X(L)$,}
\\
\hat{g} = h + dt^2 & \mbox{on $Z\times [L-1, L]$}.
\end{array}
\right.
\end{equation}
\end{Theorem}
{\bf \ref{outline}.2. Outline of the proof of Theorem \ref{main-new}.}
From now on, for simplicity we assume that $M=\emptyset$, that is $\p
W= Z$ (see Fig. \ref{cobordism}.1). The proof of the case
$M\neq\emptyset$ is rather similar to the one given below.
To make our first steps we observe the following.
\vspace{2mm}

\hspace*{20mm}
\PSbox{cob2.pstex}{10mm}{44mm}  
\begin{picture}(0,1)
\put(85,110){{\small $Z\!\times\!\!\{0\}$}}
\put(190,115){{\small $Z\!\times\!\! [1,\infty)$}}
\put(10,90){{\small $W$}}
\put(190,20){{\small $Z\!\times\!\!\{\ell\}$}}
\put(85,-13){{\small $X(\ell)$}}
\end{picture}
\vspace{5mm}

\centerline{{\small {\bf Fig. \ref{cobordism}.1.} 
The cylindrical manifold $X$ modeled by $(Z,h)$}}
\vspace{2mm}

\noindent
{\bf Observation.} We may assume that the metric $h$ on $Z$ is a
Yamabe metric with $R_h\equiv 1$. Indeed, we did not impose any
conditions on the volume of $Z$, and if there is any psc-metric
$h^{\prime}\in [h]$, then the metrics $h$ and $h^{\prime}$ are
isotopic (and, consequently, are concordant). This follows from the
fact that the set of psc-metrics $P(C)\subset C$ is convex for any
positive conformal class $C\in {\cal C}^+(Z)$ (see \cite[Lemma
7.2]{AB}).
\vspace{2mm}

\noindent
Now we start with the conformal manifold $(W,\tilde{C})$ and construct
the cylindrical manifold $X$ modeled by $(Z,h)$ as above. In
particular, we choose a conformal class $\bar{C}\in {\cal
C}^{\cyl}(X)$, with $\bar{g}\in \bar{C}$ as it was described. The idea
is to construct a function $v_{\ell}\in C^{\infty}(X)$, so that the
conformal metric $\hat{g}_{\ell}= v_{\ell}^{{4\over n-2}} \bar{g}$ on the
cylindrical manifold $X$ satisfies the conditions
\begin{equation}\label{new4}
\{
\begin{array}{l}
R_{\hat{g}_{\ell}} > 0 \ \ \ \mbox{on $X(\ell +2)$},
\\
\hat{g}_{\ell} = h + dt^2 \ \ \ \mbox{on $Z\times [\ell+1, \ell +2]$}
\end{array}
\right.
\end{equation}
for some $\ell >>1$. We achieve this in three steps.
\vspace{2mm}

\noindent
{\bf Step 1.} We study the Yamabe operator
$$
L_{\bar{g}}  = -{4(n-1)\over n-2} \Delta_{\bar{g}}  + R_{\bar{g}}
$$
on the manifold $X(\ell)$ for $\ell\geq 1$.  Namely, we study the
linear equation
\begin{equation}\label{new1}
L_{\bar{g}} u = \lambda_{\ell}u \ \ \ \mbox{on $X(\ell)$, $u|_{\p
X(\ell)}\equiv 0$,}
\end{equation}
with the Dirichet boundary condition, where $\lambda_{\ell}$ is the
corresponding first eigenvalue of $L_{\bar{g}}$ (for the Dirichet
boundary problem). For each $\ell \geq 1$ we find a function
$u_{\ell}$ satisfying (\ref{new1}) and the conditions
$$
u_{\ell}> 0 \ \ \mbox{on the interior of $X(\ell)$, and with} \ \ \ \ \ 
\min_{Z\times \{1\}} u_{\ell}= 1.
$$
In order to control the first eigenvalue $\lambda_{\ell}$, we define
the invariant $\nu_1=\nu_1(\bar{g})$ (see formula (\ref{new2}) below)
which is not a conformal invariant.  However we note (Claim \ref{l15})
that the positivity of $Y_{\bar{C}}(X(1),Z;C)$ implies positivity of
$\nu_1(\tilde{g})$ for any metric $\tilde{g}\in \bar{C}^0$. Then we
change conformally the metric $\bar{g}$ within the interior of $X(1)$
to achieve the bound $\nu_1(\bar{g})\geq 1$.
\vspace{2mm}

\noindent
{\bf Step 2.} We show that for the resulting metric $\bar{g}$ the
eigenvalues $\lambda_{\ell}$ are bounded from below:
$\lambda_{\ell}\geq \nu_1\geq 1$. Then we prove several estimates
(Claims \ref{l20}, \ref{l21}, \ref{lnew}, \ref{l25}) on the
eigenfunction $u_{\ell}$. It is important that these estimates are
independent of $\ell$.
\vspace{2mm}

\noindent
{\bf Step 3.} Here we choose a cut-off function $\phi_{\ell}$ to
define $v_{\ell}= (1-\phi_{\ell})u_{\ell+1} + \phi_{\ell}$ and examine
the scalar curvature of the conformal metric $\hat{g}$. We show that
$\hat{g}$, indeed, satisfies the conditions (\ref{new4}) for some
$\ell>>1$.

\section{Proof of Theorem \ref{main-new}}\label{cobordism}
{\bf Step 1.} We define the invariant $\nu_1=\nu_1(\bar{g})$ as follows 
\begin{equation}\label{new2}
\nu_1=\nu_1(\bar{g}) = \inf_{\begin{array}{c} _{f\in C^{\infty}(X(1))} \\
_{f\not\equiv 0}\end{array}} {\int_{X(1)} \[ {4(n-1)\over n-2}
|df|^2_{\bar{g}} + R_{\bar{g}}f^2\] d\sigma_{\bar{g}}\over
\int_{X(1)} f^2 d\sigma_{\bar{g}}}.
\end{equation}
One observes the following implication.
\begin{Claim}\label{l15}
If $Y_{\bar{C}|_{X(1)}}(X(1), Z; C) > 0 $ then $\nu_1(\tilde{g})>0$
for any metric $\tilde{g}\in \bar{C}^0$.
\end{Claim}
Thus the condition $\nu_1(\tilde{g})>0$ is conformally invariant.
\begin{Claim}\label{l16}
Under above conditions there exists a metric $\check{g}\in \bar{C}|_{X(1)}$
satisfying
$$
\{
\begin{array}{l}
\check{g} \equiv \bar{g} \ \ \mbox{on $Z\times [1-\epsilon, 1]$}
\\
\nu_1(\check{g})\geq 1
\end{array}
\right.
$$
for some $\epsilon >0$.
\end{Claim}
\begin{Proof}
We choose $\check{g}\in \bar{C}|_{X(1)}$ (keeping the condition
$\check{g}= d+ dt^2$ near $Z\times \{1\}$) so that $X(1)$ has a small
volume $\delta=\Vol_{\check{g}}(X(1))$ (see Fig. \ref{cobordism}.2).
By H\"older inequality, we have
$$
\int_{X(1)} f^2 d\sigma_{\check{g}} \leq \Vol_{\check{g}}(X(1))^{{2\over n}}
\cdot \( \int_{X(1)} |f|^{{2n\over n-2}} d\sigma_{\check{g}}\)^{{n-2\over n}}
$$
for any $f\in C^{\infty}(X(1))$. This implies
$$
{1\over \int_{X(1)} f^2 d\sigma_{\check{g}} }\geq {1\over
\delta^{{2\over n}}} \cdot {1\over  \( \int_{X(1)} |f|^{{2n\over n-2}} d\sigma_{\check{g}}\)^{{n-2\over n}}}, \ \ \mbox{and then}
$$
$$
{\int_{X(1)} \[ {4(n-1)\over n-2} |df|^2_{\check{g}} + R_{\check{g}}f^2\]
d\sigma_{\check{g}} \over \int_{X(1)} f^2 d\sigma_{\check{g}}} \geq
{1\over\delta^{{2\over n}}} \cdot {\int_{X(1)} \[ {4(n-1)\over n-2} 
|df|^2_{\check{g}}+
R_{\check{g}}f^2\] d\sigma_{\check{g}} \over \( \int_{X(1)}
|f|^{{2n\over n-2}} d\sigma_{\check{g}}\)^{{n-2\over n}}}.
$$
Thus we obtain that 
$$
\nu_1(\check{g}) \geq {1\over
\delta^{{2\over n}}} \cdot Y_{\bar{C}|_{X(1)}}(X(1), Z; [h]) ,
$$
where $Y_{\bar{C}|_{X(1)}}(X(1), Z; [h])$, perhaps, is a conformal
invariant. Finally we choose $\delta$ small enough to complete the
proof of Claim \ref{l16}.
\end{Proof}

\hspace*{3mm}
\PSbox{cob6.pstex}{10mm}{36mm}  
\begin{picture}(0,1)
\put(230,70){{\small $(W,\check{g})$}}
\put(10,70){{\small $(W,\bar{g})$}}
\end{picture}
\vspace{2mm}

\centerline{{\small {\bf Fig. \ref{cobordism}.2.} Metric $\check{g}$ on $X(1)$.}}
\vspace{1mm}

\noindent
For simplicity, we denote $\check{g}$ by $\bar{g}$. We summarize the
properties of the metric $\bar{g}\in \bar{C}$ on $X= W\cup_Z\(Z\times
[0,\infty)\)$:
\begin{equation}\label{3.1}
\{
\begin{array}{l}
\nu_1(\bar{g})\geq 1 \equiv R_h,
\\
\bar{g} = h + dt^2 \ \ \mbox{on $Z\times [1-\epsilon,\infty)$.}
\end{array}
\right.
\end{equation}
{\bf Step 2.} For any $\ell\geq 1$ we define
$$
\lambda_{\ell}= \lambda_{\ell}(\bar{g}) =
\inf_{\begin{array}{c} _{f\in C^{\infty}(X(\ell))} \\
_{f|_{\p X(\ell)}\equiv 0, \ \ f\not\equiv 0}\end{array}} 
{\int_{X(\ell)} \[ {4(n-1)\over n-2} |df|^2_{\bar{g}} +
R_{\bar{g}}f^2\] d\sigma_{\bar{g}}\over
\int_{X(\ell)} f^2 d\sigma_{\bar{g}}}.
$$
One easily proves the following  statement.
\begin{Claim}\label{l17}
The numbers $\lambda_{\ell}$ satisfy
$$
\{
\begin{array}{l}
\lambda_{\ell} \geq \min\{\nu_1, R_h\} \ \ 
\mbox{{\rm (}where $R_h\equiv 1${\rm )} for $\ell\geq 1$,}
\\
\lambda_1\geq \lambda_{\ell_1} \geq \lambda_{\ell_2} \geq 1
 \ \ 
\mbox{for $1 <\ell_1 < \ell_2$.}
\end{array}
\right.
$$
\end{Claim}
In particular, we have
\begin{equation}\label{l19}
\lambda_1\geq \lambda_{\ell} \geq R_h \equiv 1 \ \ \ \mbox{for $\ell \geq 1$}. 
\end{equation}
Now let $\ell>1$. Then there exists a function $u_{\ell}\in
C^{\infty}(X(\ell))$ such that
$$
\{
\begin{array}{l}
\displaystyle
L_{\bar{g}}u_{\ell}= - {4(n-1)\over n-2} \Delta_{\bar{g}} u_{\ell} + 
R_{\bar{g}}u_{\ell} = \lambda_{\ell} u_{\ell} \ \ \ \mbox{on $X(\ell)$},
\\
\\
u_{\ell}> 0 \ \ \ \ \ \ \ \ \ \mbox{on the interior of $X(\ell)$},
\\
\\
\displaystyle
u_{\ell}|_{\p X(\ell)} \equiv 0  \ \ \mbox{(Dirichlet boundary condition)},
\\
\\
\displaystyle
\min_{Z\times \{1\}} u_{\ell} = 1  \ \ \mbox{(normalization condition).}
\end{array}
\right.
$$
We define a function $\psi_{\ell}\in
C^{\infty}(Z\times[1,\ell])$ by
$\displaystyle
\psi_{\ell}(z,t) = 1- {t-1\over \ell-1}
$ (see Fig. \ref{cobordism}.3).

\hspace*{33mm}
\PSbox{cob7.pstex}{10mm}{25mm}  
\begin{picture}(0,1)
\put(170,15){{\small $Z\times [1,\ell]$}}
\put(17,-10){{\small $1$}}
\put(157,-10){{\small $\ell$}}
\end{picture}
\vspace{4mm}

\centerline{{\small {\bf Fig. \ref{cobordism}.3.} Function $\psi_{\ell}$}}
\vspace{2mm}

\noindent
Clearly $\Delta_{\bar{g}} \psi_{\ell} = \p_t^2
\psi_{\ell}\equiv 0$ on $Z\times [1,\ell]$. We observe the following fact.
\begin{Claim}\label{l20}
The function $u_{\ell}$ satisfies the inequality
$u_{\ell}\geq\psi_{\ell}$ on $Z\times [1,\ell]$ for any $\ell>1$.
\end{Claim}
\begin{Proof}
The above condition on $u_{\ell}$ gives that $L_{\bar{g}} u_{\ell} =
\lambda_{\ell} u_{\ell}$. This implies
\begin{equation}\label{neww1}
\begin{array}{rcl}
\Delta_{\bar{g}} u_{\ell} & =& \displaystyle
{n-2\over 4(n-1)} (R_{\bar{g}} -\lambda_{\ell}) u_{\ell} 
\\
\\
 & =& \displaystyle
{n-2\over 4(n-1)} (R_{h} -\lambda_{\ell}) u_{\ell} \leq 0 \ \ \ 
\mbox{on $Z\times [1,\ell]$}
\end{array}
\end{equation}
since $u_{\ell}\geq 0$ and by (\ref{3.1}). We obtain that
$\Delta_{\bar{g}}(\psi_{\ell}- u_{\ell}) \geq 0$. Thus by the maximum
principle,
$$
\psi_{\ell}- u_{\ell} \leq \max_{\p( Z\times [1,\ell])} (\psi_{\ell}-
u_{\ell}) = 0
$$
because of the choice of $\psi_{\ell}$ (here, of course, $\p( Z\times
[1,\ell]) = (Z\times \{1\})\sqcup (Z\times \{\ell\})$). Then we obtain that
$u_{\ell}\geq\psi_{\ell}$ on $Z\times [1,\ell]$.
\end{Proof}
\begin{Claim}\label{l21}
Let $\ell >3$. Then there exists a constant $K>0$ independent of
$\ell$, so that
$$
|d u_{\ell}| \leq K \cdot |u_{\ell}| \ \ \ \mbox{on $Z\times [1,\ell-2]$.}
$$
\end{Claim}
\begin{Proof}
Let $x\in X$. We define $ B_1(x) = \{ y\in X \ | \
\dist_{\bar{g}}(y,x)\leq 1 \ \}.  $ Now we have to recall the
following facts.
\begin{Fact}\label{l22}
{\rm (see \cite{Jost})}\\ Let $f\in
C^{\infty}(X)$. Then for any $x\in X$ there exists a constant $K_1>0$
so that
$$
|df(x)| \leq K_1\[ \int_{ B_1(x)} |\Delta_{\bar{g}} f| \ 
d\sigma_{\bar{g}} + \int_{ B_1(x)} |f| \ 
d\sigma_{\bar{g}} + \int_{ \p B_1(x)} |f| \ 
d\sigma_{\bar{g}|_{ \p B_1(x)}}\].
$$
\end{Fact}
\begin{Fact}\label{l23}
{\rm (Harnack inequality, see \cite[Theorem 8.20]{GT})}\\ There exists
a constant $K_2>0$ independent of $\ell>1$ so that
$$
|u_{\ell}|\leq K_2 \cdot \( \inf_{y\in B_1(x)} u_{\ell}(y)\)
\ \ \ \ 
\mbox{for $x\in Z\times [1,\ell-1]$.} 
$$
\end{Fact}
We continue with the proof of Claim \ref{l21}. 
The Facts \ref{l22}, \ref{l23}, (\ref{l19}) and (\ref{neww1}) imply
that there exists a constant  $K_3>0$ so that
$$
\begin{array}{rcl}
|d u_{\ell}(x)| & \leq & \displaystyle
K_1 \cdot \[ K_3 \int_{ B_1(x)} |u_{\ell}| 
d\sigma_{\bar{g}} + \int_{ B_1(x)} |u_{\ell}| 
d\sigma_{\bar{g}}  + \int_{ \p B_1(x)} |u_{\ell}| \ 
d\sigma_{\bar{g}|_{ \p B_1(x)}}\]
\\
\\
& \leq & \displaystyle
K_1K_2(K_3+2) \cdot   \( \inf_{y\in B_1(x)} u_{\ell}(y)\) 
\end{array}
$$
for $x\in Z\times [1,\ell-2]$. This completes the proof of Claim \ref{l21}.
\end{Proof}
Recall that $\displaystyle \min_{Z\times \{1\}} u_{\ell} = 1 $. Now
Claim \ref{l21} implies that there exists a constant $\bar{K}>0$,
independent of $\ell$, so that
$\displaystyle 
\max_{Z\times \{1\}} u_{\ell} \leq \bar{K}.
$
\begin{Claim}\label{l25}
The function $u_{\ell}$ satisfies 
$$
u_{\ell} \leq \bar{K} \ \ \ \mbox{on $Z\times [1, \ell]$ for any $\ell>1$.}
$$
\end{Claim}
\begin{Proof}
Consider the function $e^{\delta t}\cdot u_{\ell}$ for $\delta
>0$. Recall that
$$
\Delta_{\bar{g}} u_{\ell} = {n-2\over 4(n-1)} (R_h - \lambda_{\ell}).
$$
We use (\ref{l19}) and Claim \ref{l21} to see the following estimate:
$$
\begin{array}{rcl}
\Delta_{\bar{g}}(e^{\delta t}\cdot u_{\ell}) & = & \displaystyle 
e^{\delta t} \cdot \Delta_{\bar{g}} u_{\ell} + 
2\delta \cdot e^{\delta t}  u_{\ell}^{\prime} +  
\delta^2\cdot e^{\delta t}   u_{\ell}
\\
\\
& \geq & \displaystyle 
-{n-2\over 4(n-1)} \lambda_1 \cdot e^{\delta t}  u_{\ell} - 2 \delta K
e^{\delta t}  u_{\ell} + \delta^2\cdot e^{\delta t}   u_{\ell}
\\ 
\\ 
& \geq & \displaystyle \(\delta^2 -{n-2\over 4(n-1)}\lambda_1 -
2 \delta K\) e^{\delta t} u_{\ell} \geq 0
\end{array}
$$
for large enough $\delta >> 1$ (since $u_{\ell}\geq 0$), where
$(\cdot)^{\prime} = {\p \over \p t}(\cdot)$.  Now the maximum
principle gives
$$
e^{\delta t} u_{\ell} \leq \max_{\p (Z\times [1,\ell])}  e^{\delta t}
u_{\ell} {=} \max_{Z\times \{ 1 \}}  e^{\delta t} u_{\ell}
\leq e^{\delta t} \cdot \bar{K}
$$
on $Z\times [1,\ell]$ since $u_{\ell}|_{Z\times \{\ell\}}=0$.  Thus we
obtain $ u_{\ell} (z,t) \leq e^{\delta(1- t)}\cdot \bar{K} \leq
\bar{K} $ on $Z\times [1,\ell]$.
\end{Proof}
We need one more precise estimate on the function $u_{\ell}$.
\begin{Claim}\label{lnew}
There exist constants $\tilde{K}_1>0$, $\tilde{K}_2>0$ (independent of
$\ell$) such that
$$
|d u_{\ell}|_{C^{0,\alpha}} \leq \tilde{K}_1\( |u_{\ell}| + \tilde{K}_2\)
\ \ \mbox{on $Z\times [1,\ell]$}.
$$
\end{Claim}
\begin{Proof}
Indeed, we have that the function $u_{\ell}$ satisfies
$$
\{
\begin{array}{l}
u_{\ell}|_{Z\times \{\ell\]} \equiv 0,
\\
0\leq u_{\ell} \leq \bar{K} \ \ \mbox{on $Z\times [1,\ell]$},
\\
|d u_{\ell}| \leq K\cdot |u_{\ell}| \ \ \mbox{on $Z\times [1,\ell-2]$}.
\end{array}
\right.
$$
Recall that we have
$$
\begin{array}{l}\displaystyle
\Delta_{\bar{g}} u_{\ell} = {n-2\over 4(n-1)}(R_{\bar{g}} -\lambda_{\ell})
\ \ \mbox{on $Z\times [1,\ell]$}
\\
\\
\displaystyle
\lambda_1 \geq \lambda_{\ell} \geq R_h \equiv 1.
\end{array}
$$
Then, by standard argument, we obtain that
$$
|d u_{\ell}|_{C^{0,\alpha}} \leq \tilde{K}\cdot |u_{\ell}| 
\ \ \mbox{on $Z\times [1,\ell-2]$}.
$$
Then \cite[Theorem 8.33]{GT} implies that there exist
constants $\tilde{K}_1>0$, $\tilde{K}_2>0$  such that
$$
|d u_{\ell}|_{C^{0,\alpha}} \leq \tilde{K}_1\( |u_{\ell}| + \tilde{K}_2\)
\ \ \mbox{on $Z\times [1,\ell]$}.
$$
This completes the proof of Claim \ref{lnew}.
\end{Proof}

\hspace*{33mm}
\PSbox{cob8.pstex}{10mm}{22mm}  
\begin{picture}(0,1)
\put(200,15){{\small $Z\times [0,\ell]$}}
\put(17,-10){{\small $1$}}
\put(162,-10){{\small $\ell+1$}}
\put(85,-10){{\small $1+ {\ell\over 2}$}}
\put(-40,50){{\small $1$}}
\put(-40,25){{\small ${1\over 2}$}}
\end{picture}
\vspace{4mm}

\centerline{{\small {\bf Fig. \ref{cobordism}.4.} Function $\phi_{\ell}$}}
\vspace{2mm}

\noindent
{\bf Step 3.} Let $\ell>> 1$. Let $\phi_{\ell}\in C^{\infty}(X)$ be a
cut-off function satisfying the following conditions (see
Fig. \ref{cobordism}.4):
\begin{enumerate}
\item[{\bf (1)}]
$
\phi_{\ell}(x) = \{
\begin{array}{ll}
0 & \mbox{for $x\in X(1)$},
\\
1 & \mbox{for $x\in Z\times [\ell+ 1,\infty)$}
\end{array}
\right.
$
\item[{\bf (2)}]
$0 \leq \phi_{\ell} \leq 1$ on $X$, and $\phi_{\ell}(z,t)= \phi_{\ell}(t)$
for $(z,t)\in Z\times [1, \ell+ 1]$.
\item[{\bf (3)}] For some constant $\hat{K}>0$ (independent of $\ell$)
$$
0\leq \phi^{\prime}_{\ell} \leq {\hat{K}\over \ell},
\ \ \ \ | \phi^{\prime\prime}_{\ell} | \leq {\hat{K}\over \ell^2}
\ \ \ \ 
\mbox{on $Z\times [1,\ell+1]$}.
$$
\item[{\bf (4)}] Moreover, $\phi_{\ell}(1+{\ell\over 2})= {1\over 2}$ on
$Z\times \{ 1+{\ell\over 2}\}$.
\end{enumerate}
\noindent
It is not difficult to find such function $\phi_{\ell}$. We let
$v_{\ell}= (1-\phi_{\ell})\cdot u_{\ell+1} + \phi_{\ell}\in
C^{\infty}_+(X)$, and the conformal metric $
\hat{g}_{\ell}= v_{\ell}^{{4\over n-2}}\cdot \bar{g} $ on $X$. Then
the scalar curvature of the metric $\hat{g}_{\ell}$ is given by
$$
R_{\hat{g}_{\ell}} = v_{\ell}^{-{n+2\over n-2}} \[ -{4(n-1)\over n-2}
\Delta_{\bar{g}} v_{\ell} + R_{\bar{g}} v_{\ell}\] = 
v_{\ell}^{-{n+2\over n-2}} L_{\bar{g}} v_{\ell}.
$$
We examine the scalar curvature $R_{\hat{g}_{\ell}}$
on three different pieces: 
$$
X = X(1)\cup \( Z\times [1, \ell+1]\) \cup \( Z\times [\ell+1, \infty)\).
$$
$\bullet$ The piece $X(1)$.  Then we have that $v_{\ell} \equiv
u_{\ell+1}$, thus
$$
R_{\hat{g}_{\ell}}|_{X(1)}= u_{\ell+1}^{-{n+2\over n-2}}\cdot
(\lambda_{\ell+1}u_{\ell+1}) = 
\lambda_{\ell+1} \cdot u_{\ell+1}^{-{4\over n-2}} > 0 .
$$
$\bullet$ The piece $Z\times [\ell+1, \infty)$. Here we have that
$v_{\ell}\equiv 1$ (which is equivalent to the fact that
$\hat{g}_{\ell}= \bar{g} = h+ dt^2$). Thus
$$
R_{\hat{g}_{\ell}}|_{Z\times [\ell+1, \infty)} = R_{\bar{g}} = R_h \equiv 1.
$$
$\bullet$ The piece $Z\times [1, \ell+1]$. This case is more
complicated. We have:
$$
\begin{array}{rcl}
R_{\hat{g}_{\ell}}|_{Z\times [1, \ell+1]} &=&\displaystyle
v_{\ell}^{-{n+2\over n-2}} \[ (1-\phi_{\ell}) L_{\bar{g}} u_{\ell+1} +
{4(n-1)\over n-2}(u_{\ell+1} -1)\phi^{\prime\prime}_{\ell}\right.
\\
\\
& &\displaystyle
\left. \ \ \ \ \ \ \ \ \ \ \ \ \ \ \ \ \ \ \ \ \ \ \ \ \ \ 
\ \ \ \ \ \ \ \ \ \ \ \ \ \ \ \ \ \ \ \ \ \ \ \ \ \ 
 + {8(n-1)\over n-2} \phi^{\prime}_{\ell}\cdot u_{\ell+1}^{\prime} + 
R_h \phi_{\ell}\]
\\
\\
&=&\displaystyle
v_{\ell}^{-{n+2\over n-2}} \[ (1-\phi_{\ell}) \lambda_{\ell+1}\cdot u_{\ell+1}+
{4(n-1)\over n-2}(u_{\ell+1} -1)\phi^{\prime\prime}_{\ell}\right.
\\
\\
& &\displaystyle
\left. \ \ \ \ \ \ \ \ \ \ \ \ \ \ \ \ \ \ \ \ \ \ \ \ \ \ 
\ \ \ \ \ \ \ \ \ \ \ \ \ \ \ \ \ \ \ \ \ \ \ \ \ \ 
 + {8(n-1)\over n-2} \phi^{\prime}_{\ell}\cdot u_{\ell+1}^{\prime} + 
\phi_{\ell}\].
\end{array}
$$
We use (\ref{l19}),   Claim \ref{l25}, Claim \ref{lnew}, 
and property (3) of the function $\phi_{\ell}$ to get the estimation:
$$
R_{\hat{g}_{\ell}}|_{Z\times [1, \ell+1]} \geq v_{\ell}^{-{n+2\over n-2}} 
\[
(1-\phi_{\ell}) u_{\ell+1} -{8\bar{K}\hat{K}\over \ell^2} -
{8\bar{K}\tilde{K}_1( \hat{K} + \tilde{K}_2)\over \ell} + \phi_{\ell}
\].
$$
Now we examine even more carefully the scalar curvature
$R_{\hat{g}_{\ell}}$ on the cylinder
$$
\begin{array}{c}
Z\times [1, \ell+1] = \(Z\times [1, 1+ {\ell\over 2}]\)\cup 
\(Z\times [1+ {\ell\over 2}, \ell+1]\). 
\end{array}
$$
$\bullet$ The piece $Z\times [1, 1+ {\ell\over 2}]$. Here the property
(4) of the function $\phi_{\ell}$ and Claim \ref{l20} imply:
$$
\{
\begin{array}{l}
1-\phi_{\ell}\geq {1\over 2},
\\
u_{\ell+1}\geq \psi_{\ell+1} \geq \psi_{\ell+1} (1+{\ell\over 2})
\geq {1\over 2}.
\end{array}
\right.
$$
Thus we have that
$$
R_{\hat{g}_{\ell}}|_{Z\times [1, 1+{\ell\over 2}]} \geq
v_{\ell}^{-{n+2\over n-2}} \({1\over 4} - {8\bar{K}\hat{K}\over \ell^2} - 
{8\bar{K}\tilde{K}_1( \hat{K} + \tilde{K}_2)\over \ell} \).
$$
Clearly there exists such $\ell_1>>1$ that $R_{\hat{g}_{\ell}}>0$ on
$Z\times [1, 1+ {\ell\over 2}]$ for all $\ell\geq \ell_1$.
\vspace{2mm}

\noindent
$\bullet$ The piece $Z\times [1+ {\ell\over 2},\ell+1]$. Here we have
$\phi_{\ell}\geq {1\over 2}$ by the conditions (3) and (4) on the
function $\phi_{\ell}$. Thus we we have
$$
R_{\hat{g}_{\ell}}|_{Z\times [1+ {\ell\over 2},\ell+1]} \geq
v_{\ell}^{-{n+2\over n-2}} \({1\over 2} - {8\bar{K}\hat{K}\over \ell^2} - 
{8\bar{K}\tilde{K}_1( \hat{K} + \tilde{K}_2)\over \ell}\).
$$
Thus there exists $\ell_2>>1$ such that $R_{\hat{g}_{\ell}}> 0$
on $Z\times [1+ {\ell\over 2},\ell+1]$ for all $\ell\geq\ell_2$.
\vspace{2mm}

\noindent
Now let $\ell_0 = \max\{\ell_1,\ell_2\}$, and let $\hat{g}=
\hat{g}_{\ell_0}$, and $L= \ell_0+2$. Thus we constructed a metric
$\hat{g} \in \Riem(X(L))$ such that
$$
\{
\begin{array}{ll}
R_{\hat{g}}> 0 & \mbox{on $X(L)$},
\\
\hat{g} = h+ dt^2 & \mbox{on $Z\times [L-1, L]$}.
\end{array}
\right.
$$
This completes the proof of Theorem \ref{main-new}. \hfill $\Box$
\section{Some topology}\label{topology}
{\bf \ref{topology}.1. Summary on $\gamma$-structures.} We briefly
review necessary definitions and constructions given by S. Stolz
\cite{Stolz1}. Let $\pi=\pi_1(M)$ be the fundamental group of $M$, and 
$w_i(M)\in H^i(M;\Z_2)$ be the Stiefel-Whitney characteristic classes.
\vspace{2mm}

\noindent
The main conceptual issue here is to determine precisely which
topological structure on a smooth compact manifold $M$ carries
complete information on the existence of a psc-metric on $M$. Indeed,
it is well-known that the fundamental group $\pi$ is crucially
important for the existence question. Then there is clear difference
when a manifold $M$ is oriented or not (which depends on $w_1(M)$). On
the other hand, a presence of the $Spin$-structure (which means that
$w_2(M)=0$) gives a way to use the Dirac operator on $M$ to control
the scalar curvature via the vanishing formulas. S. Stolz puts
together those invariants to define a $\gamma$-structure.
\vspace{2mm}

\noindent
To simplify our presentation, we consider only the case of oriented
manifolds.  In the oriented case the $\gamma$-structures have very
transparent geometric description. The non-oriented case is more
subtle and complicated; we would like to live this case outside of our
paper.  However, in our view, one should not meet any difficulties to
generalize our results to the non-oriented case.
\vspace{2mm}

\noindent
Let $M$ be an oriented manifold with $\pi=\pi_1(M)$. Let $f : M \lra
B\pi$ be a classifying map for the fundamental group, and $p:
\widetilde{M} \rightarrow M$ be the universal cover. Recall that the
second Stiefel-Whitney class $w_2=w_2(M)$ is zero if and only if the
manifold $M$ admits a $Spin$ structure. We have the following three
cases to consider:
\begin{enumerate}
\item[{\bf (1)}] $w_2=0$, thus the manifold $M$ is a $Spin$
manifold;
\item[{\bf (2)}] $w_2\neq 0$, but the universal cover
$\widetilde{M}$ is a $Spin$ manifold;
\item[{\bf (3)}] $w_2\neq 0$, and the universal cover
$\widetilde{M}$ is not a $Spin$ manifold.
\end{enumerate}
{\bf Comments. (1)} In this case $M$ admits a $Spin$ structure,
however it is important to choose the $Spin$ structure. We call a
manifold $M$ a $Spin$-manifold if the $Spin$-structure is chosen.  A
classifying map $f : M \lra B\pi$ then determines a canonical
cobordism class $[(M,f)]\in \Omega^{\Spin}_n(B\pi)$. 
In this case a $\gamma$-structure on $M$  is defined as a choosen
$Spin$-structure on $M$ together with the classifying map $f : M \lra
B\pi$ for the fundamental group $\pi$. 
\vspace{2mm}

\noindent
{\bf (2)} This case involves  more. Consider the induced
homomorphism 
$$
f^*: H^2(B\pi;\Z_2) \rightarrow H^2(M;\Z_2).
$$ 
In this case S. Stolz proves that there exists a unique element $e\in
H^2(B\pi;\Z_2)$, so that $f^*(e) = w_2$. The element $e$, as any
element of $H^2(B\pi;\Z_2)$, determines a central group extension $1
\rightarrow \Z_2 \rightarrow \hat{\pi} \stackrel{\rho}{\rightarrow}
\pi \rightarrow 1$. Futhermore, this extension splits (or trivial) if
and only if $e=0$. Thus the pair $(\hat{\pi},\pi)$ completely encodes
the case (2) and the case (1) as well (then $e=0$, and $\hat{\pi}\cong
\pi\times \Z_2$). This gives the $\gamma$-structure
$\gamma=(\hat{\pi},0,\pi)$ in the notations of
\cite{Stolz1}. Alternatively, this structure gives the following
construction. Let $\sigma\in \Z_2 \subset \hat{\pi}$ be a generator.
Then the element $(\sigma,-1)$ is central in the direct product
$\hat{\pi}\times Spin(n)$. The Lie group $G(\gamma,n)$ is defined as a
factor group of $\hat{\pi}\times Spin(n)$ by the central subgroup
$\Z_2$ (generated by $(\sigma,-1)$). By construction, the group
$G(\gamma,n)$ has a canonical homomorphism $j : G(\gamma,n)\lra
SO(n)$.
\vspace{2mm}

\noindent
Now let $g$ be a Riemannian metric on $M$, then a chosen orientation
on $M$ gives the frame bundle $P_{SO(n)}(M) \rightarrow M$. S. Stolz
shows \cite{Stolz1} that in this case the $\gamma$-structure
determines a canonical principal bundle $P_{G(\gamma,n)}(M) \lra
P_{SO(n)}(M)$. We obtain the principal bundle $P_{G(\gamma,n)}(M) \lra
M$, and thus a map $ \hat{f} : M \lra BG(\gamma,n) $ to the
classifying space. The case when the above extension $e$ is trivial
gives the isomorphism $G(\gamma,n)\cong \pi\times Spin(n)$. Otherwise
the group $G(\gamma,n)$ is a ``twisted (by the extension $e$)
version'' of the group $Spin(n)$.
\vspace{2mm}

\noindent
We remark that the group $G(\gamma,n)$ determines the Thom space
$MG(\gamma,n)$, and thus the cobordism groups $\Omega_n(\gamma)$ given
via the Thom-Pontryagin construction. In particular, the pair
$(M,\hat{f})$ determines a cobordism class in $\Omega_n(\gamma)$
(where $n=\dim M$). Both cases (1) and (2) are described in
\cite{Stolz1} as $\gamma=(\hat{\pi}, 0,\pi)$ with $\hat{\pi}$ given by
the above extension $e$.
\vspace{2mm}

\noindent
{\bf (3)} This case is easy. The $\gamma$-structure here is nothing
but a choice of orientation on $M$ together with the classifying map
$f : M \lra B\pi$. Then the pair $(M,f)$ gives a cobordism class in
the oriented cobordism ring $\Omega_n^{\SO}(B\pi)$. In the notations
of
\cite{Stolz1}, $\gamma=(\pi,0,\pi)$. 
\vspace{2mm}

\noindent
{\bf Conclusion.} We emphasize that in each of the above cases we have
a well-defined cobordism category ${\mathfrak M}(\gamma)$ of manifolds
equipped with $\gamma$-structure. Let $M_0$, $M_1$ be two manifolds
equipped with given $\gamma$-structure. A cobordism $W$ between $M_0$
and $M_1$ in the category ${\mathfrak M}(\gamma)$ is called
$\gamma$-cobordism. 
\vspace{2mm}

\noindent
{\bf \ref{topology}.2. Conformal and psc-cobordism groups.} Now let
$\PPos(\gamma)$ be the following cobordism category. The objects of
$\PPos(\gamma)$ are the pairs $(M,g)$, where $M$ is a manifold with
$\gamma$-structure, and $g$ is a psc-metric on $M$. Manifolds
$(M_0,g_0)$, $(M_1,g_1)$ are psc-cobordant in the category
$\PPos(\gamma)$ if they there is a $\gamma$-cobordism $W$ between $M_0$
and $M_1$, where $W$ is given a psc-metric $\bar{g}$, so that
$$
R_{\bar{g}}>0, \ \ \ \ 
\bar{g}|_{M_j}=g_j, \ \ \ \ \mbox{and $\bar{g}=g_j +dt^2$ near the
boundary $\p W= M_0\sqcup (-M_1)$ for $j=0,1$.}
$$
We denote the corresponding cobordism groups $\Pos_n(\gamma)$.  We
emphasize that we restrict our attention to the dimensions $n\geq 5$.
\vspace{2mm}

\noindent
The corresponding conformal cobordism category $\PPos^{\conf}(\gamma)$
is defined similarly. The objects of $\PPos^{\conf}(\gamma)$ are
positive conformal $\gamma$-manifolds $(M,C)$, where, as before, $M$
is a manifold with $\gamma$-structure, and $C\in {\cal C}^+(M)$ is a
positive conformal class. Then two positive conformal
$\gamma$-manifolds $(M_0,C_0)$, $(M_1,C_1)$ are conformally cobordant
if there exists a $\gamma$-cobordism $W$ between $M_0$ and $M_1$, so
that the relative Yamabe invariant $ Y (W,M_0\sqcup (-M_1); C_0\sqcup
C_1)>0$. 
\vspace{2mm}

\noindent
We denote the corresponding cobordism groups
$\Pos_n^{\conf}(\gamma)$. Here we also let $n\geq 5$ (however, all
definitions make sense for $n=2,3,4$ as well). The fact, that the
conformal cobordism is an equivalence relation is not entirely trivial
(see proof in \cite{AB}). A group structure here is given by taking a
disjoint union of manifolds.  We have a canonical functor
$\PPos(\gamma) \lra \PPos^{\conf}(\gamma)$ given by taking conformal
classes of corresponding metrics, so we have natural homomorphism
$\Pos_n(\gamma)\lra \Pos_n^{\conf}(\gamma)$. Clearly Theorem
\ref{main-new} implies Theorem A and Corollary B. Since concordance is
just a particular case of cobordism, this also implies Corollary D.
\vspace{2mm}

\noindent
{\bf \ref{topology}.3. Relative cobordism groups.}  Now we define the
cobordism category ${\mathfrak R}(\gamma)$ for a given
$\gamma$-structure as above. The objects of the category ${\mathfrak
R}(\gamma)$ are $\gamma$-manifolds $(M,\p M; \bar{g}, g)$, where
$\bar{g}$ is a Riemannian metric on $M$, and $g$ is a psc-metric on
$\p M$, such that
$$
\mbox{$\bar{g}=g +dt^2$ near the boundary $\p M$.}
$$
In particular, if $M$ is a closed $\gamma$-manifold, and $\bar{g}$ is
any Riemannian metric, then $(M,\emptyset;\bar{g},\emptyset)$ is an
object of ${\mathfrak R}(\gamma)$.  Two manifolds $(M_0,\p
M_0;\bar{g}_0,g_0)$, $(M_1,\p M_1;\bar{g}_1, g_1)$ like this are
cobordant in the category ${\mathfrak R}(\gamma)$ if there exist a
$\gamma$-manifold $(W,\p W; \tilde{g}, \hat{g})$ with given
decomposition of the boundary
$$
\p W = M_0 \cup_{\p M_0} V\cup_{\p M_1} (- M_1), \ \ \ 
$$
where $\p V = \p M_0 \sqcup (-\p M_1)$, such that (see
Fig. \ref{topology}.1)
$$
\begin{array}{ll}
{\mathbf (a)} & \hat{g}|_{\p V} = g_0 \sqcup g_1, \ \ 
\mbox{with $\hat{g}= \hat{g}|_{\p V} +dt^2$ near $\p V$,}
\\
{\mathbf (b)} &\mbox{$R_{\hat{g}}>0$ on $V$},
\\
{\mathbf (c)} & \tilde{g}|_{\p W} = \hat{g} = \bar{g}_0\cup \hat{g}|_{V} 
\cup \bar{g}_1, \ \ \mbox{and}
\\
{\mathbf (d)} & \tilde{g} = \hat{g} + dt^2  \ \ 
\mbox{near the boundary $\p W$.}
\end{array} \ \ \ \ \ \ \ \ \ \ \  \ \ \ \ \ \ \ 
\ \ \ \ \ \ \ \ \ \ \  \ \ 
$$
Here ``$-M$'' means the same manifold $M$ with the choice of {\sl
opposite $\gamma$-structure} (see \cite{Stolz1} for more details).  We
remark that the manifold $(V, \hat{g}|_{V})$ delivers a psc-cobordism
between $(\p M_0, g_0)$ and $(\p M_1, g_1)$ (we emphasized this by a
bold line in Fig. \ref{topology}.1).
\vspace{2mm}

\hspace*{23mm}
\PSbox{cob9.pstex}{10mm}{50mm}  
\begin{picture}(0,1)
\put(60,40){{\small $(W,\tilde{g})$}}
\put(-35,117){{\small $(M_0,\bar{g}_0)$}}
\put(30,117){{\small $\p M_0$}}
\put(70,117){{\small $(V, \hat{g}|_V)$}}
\put(115,117){{\small $\p M_1$}}
\put(165,117){{\small $(-M_1,\bar{g}_1)$}}
\end{picture}
\vspace{2mm}

\centerline{{\small {\bf Fig. \ref{topology}.1.} Cobordism in the
category ${\mathfrak R}(\gamma)$}}
\vspace{1mm}

\noindent
Again, we emphasize that for each Riemannian manifold with boundary it
is assumed here that a metric is a product metric near its boundary.
Let $R_n(\gamma)$ be the corresponding cobordism groups. Disjoint
union of manifolds induces an abelian group structure on $R_n(\gamma)$
(see \cite{Stolz1}).
\vspace{2mm}

\noindent
The conformal cobordism category ${\mathfrak R}^{\conf}(\gamma)$ is
defined similarly. To avoid any confusions, we spell out the
definition.  The objects of ${\mathfrak R}^{\conf}(\gamma)$ are
conformal $\gamma$-manifolds $(M,\p M; \bar{C},C)$, where
$(\bar{C},C)\in {\cal C}(M,\p M)$ (i.e. $\p \bar{C}= C$
$\Longleftrightarrow$ $\bar{C}|_{\p M}= C$) with $C\in {\cal C}^+(\p
M)$ positive conformal class. Conformal manifolds $(M_0,\p M_0;
\bar{C}_0,C_0)$, $(M_1,\p M_1; \bar{C}_1,C_1)$ like this are cobordant
in the category ${\mathfrak R}^{\conf}(\gamma)$ if there is a
conformal manifold $(W,\p W, \tilde{C}, \hat{C})$ with given
decomposition of the boundary
$$
\p W = M_0 \cup_{\p M_0} V\cup_{\p M_1} (- M_1), \ \ \ 
$$
where $\p V = \p M_0 \sqcup (-\p M_1)$, such that 
$$
\begin{array}{ll}
(a)^{\conf} & \hat{C}|_{\p V} = C_0 \sqcup C_1, 
\\
(b)^{\conf} & Y_{\hat{C}|_V} (V, \p M_0 \sqcup \p M_0; C_0\sqcup C_1)>0,
\\
(c)^{\conf} & \tilde{C}|_{\p W} = \hat{C} = \bar{C}_0\cup \hat{C}|_{V} 
\cup \bar{C}_1.
\end{array} \ \ \ \ \ \ \ \ \ \ \  \ \ \ \ \ \ \ 
\ \ \ \ \ \ \ \ \ \ \  \ \ \ \ \ \ \ \ \ \ \ \ \  \ \ 
$$
Let $R_n^{\conf}(\gamma)$ be the corresponding cobordism groups.
Clearly there are natural homomorphisms $j: \Omega_n(\gamma) \lra
R_n(\gamma)$ and $j^{\prime}: \Omega_n(\gamma) \lra
R_n^{\conf}(\gamma)$, given by assigning an arbitrary Riemannian
metric (or conformal class to a $\gamma$-manifold). We remark here
that two closed $\gamma$-manifolds $(M,g_0)$ and $(M, g_1)$ (with any
two metrics $g_0$, $g_1$) are cobordant in the category ${\mathfrak
R}(\gamma)$ since the space of Reimannian metrics is convex. Thus a
linear homotopy $g_t = (1-t)g_0 + tg_1$ gives a metric on the cylinder
$M\times [0,1]$. The same is true for conformal manifolds if we do not
impose any conditions on conformal classes. The maps $\p: R_n(\gamma)
\lra Pos_{n-1}(\gamma)$ and $\p^{\prime}: R_n^{\conf}(\gamma) \lra
Pos_{n-1}^{\conf}(\gamma)$ are given by taking all data on
boundaries. Finally one has the forgeting (metric or conformal class)
homomorphisms $F: Pos_{n}(\gamma) \lra \Omega_n(\gamma)$, and
$F^{\prime}: Pos_{n}^{\conf}(\gamma) \lra \Omega_n(\gamma)$. It is
easy to show that the following diagram is commutative and has exact
rows:
$$
\begin{diagram}
\setlength{\dgARROWLENGTH}{.8em}
	\node{\cdots}
		\arrow{e}
	\node{R_{n+1}(\gamma)} 
		\arrow{e,t}{\p}
		\arrow{s,l}{c}
	\node{Pos_{n}(\gamma)} 
		\arrow{e,t}{F}
		\arrow{s,l}{\cong}
	\node{\Omega_{n}(\gamma)} 
		\arrow{e,t}{j}
		\arrow{s,l}{Id}
	\node{R_{n}(\gamma)} 
		\arrow{e,t}{\p}
		\arrow{s,l}{c}		
	\node{Pos_{n-1}(\gamma)} 
		\arrow{e}
		\arrow{s,l}{\cong}
	\node{\cdots}
\\
	\node{\cdots}
		\arrow{e}
	\node{R_{n+1}^{\conf}(\gamma)} 
		\arrow{e,t}{\p^{\prime}}
	\node{Pos_{n}^{\conf}(\gamma)} 
		\arrow{e,t}{F^{\prime}}
	\node{\Omega_{n}(\gamma)} 
		\arrow{e,t}{j^{\prime}}
	\node{R_{n}^{\conf}(\gamma)} 
		\arrow{e,t}{\p^{\prime}}	
	\node{Pos_{n-1}^{\conf}(\gamma)} 
		\arrow{e}
	\node{\cdots}
\end{diagram}
$$
Five-lemma implies that $c : R_{n}(\gamma) \lra R_{n}^{\conf}(\gamma)$
is an isomorphism. This concludes the proof of Corollary C.

\section{Discussion}\label{disscussion}
{\bf \ref{disscussion}.1. Concordance classes and groups
$R_n(\gamma)$.} To make our discussion transparent, we concentrate our
attention on the case of simply connected $Spin$ manifolds, then
$\Omega_n(\gamma)=\Omega^{\Spin}_n$.  In this case we omit the
``$\gamma$-notation'' for all cobordism groups we have
here. Futhermore, we consider the simplest possible manifold, the
standard sphere. Thus let $M = S^n$ for $n\geq 5$, and let $\Pi_n$ be
the set of psc-concordant classes of psc-metrics on $S^n$. Corollary
D, in particular, identifies the set $\Pi_n$ with its ``conformal
relative'', the set $\Pi_n^{\conf}$ of conformally concordant positive
conformal classes on $S^n$.  The connective sum operation induces an
abelian group structure on $\Pi_n$ with zero class represented by the
standard metric $g_{\can}$. Thus $\Pi_n^{\conf}$ inherits this group
structure.
\vspace{2mm}

\noindent
On the other hand, it is known (see \cite{Hajduk}) that for simply
connected $Spin$ manifolds the relative psc-cobordism groups $R_n$ are
naturally isomorphic to the concordance groups $\Pi_n$. We obtain the
isomorphisms: $ R^{\conf}_n\cong R_n \cong \Pi_n \cong \Pi_n^{\conf}$.
\vspace{2mm}

\noindent
{\bf Conlusion.} {\sl The groups $R^{\conf}_n$ (and, consequently, the
groups $\Pos_n^{\conf}$) are completely determined by the concordance
classes of a single manifold: the sphere $S^n$.} 
\vspace{2mm}

\noindent
In the conformal world, we have the set $\pi_0({\cal C}^+(S^n))$ of
positive conformal classes.
\vspace{2mm}

\noindent
{\bf Problem 1.} {\sl Study the relationship of the group of conformal
concordance classes $\Pi_n^{\conf}$ and $\pi_0({\cal C}^+(S^n))$.}
\vspace{2mm}

\noindent
A study of the space ${\cal C}^+(S^n)$ naturally leads to an
interesting model of moduli space of positive conformal classes. 
\vspace{2mm}

\noindent
{\bf \ref{disscussion}.2. Moduli spaces.} Again, we consider the
sphere $S^n$ with $n\geq 5$. A standard definition of the moduli space
of psc-metrics goes as follows.  Let $\Riem^+(S^n)\subset \Riem(S^n)$
be the space of psc-metrics. The diffeomorphism group $\Diff_+(S^n)$
of orientation-preserving diffeomorphisms of the sphere $S^n$
naturally acts (on the right) on the space of metrics $\Riem(S^n)$ by
pulling back a metric. Obviously this action preserves the subspace
$\Riem^+(S^n)$. There is a serious problem with this action: it is far
away from to be free, leaving us very little chance to understand the
topology of the moduli space ${\cal M}^{+}(S^n) =
\Riem^+(S^n)/\Diff_+(S^n)$ of psc-metrics.
\vspace{2mm}

\noindent
We would like to suggest an alternative construction of such a moduli
space following the paper \cite{Morava} by J. Morava \& H. Tamanoi.
The construction below holds for arbitrary compact smooth manifold,
not just for the sphere $S^n$.
\vspace{2mm}

\noindent
Let ${\cal C}(S^n)$, ${\cal C}^+(S^n)$ be the spaces of all conformal
classes and positive ones. The projection map $\Riem(S^n)\lra {\cal
C}(S^n)$ induces the map $\Riem^+(S^n)\lra {\cal C}^+(S^n)$. Clearly
both spaces $\Riem(S^n)$ and ${\cal C}(S^n)$ are contractible, and
again the diffeomorphism group $\Diff_+(S^n)$ action on ${\cal
C}(S^n)$ is not free. To refine the construction, we choose a base
point $x_0\in S^n$.
\vspace{2mm}

\noindent
The space of conformal classes ${\cal C}(S^n)$ is the orbit space
of the action (left multiplication) of the group $C_+^{\infty}(S^n)$ on
the space of metrics $\Riem(S^n)$. With a given base point $x_0\in S^n$,
we consider the following subspace of
$C_+^{\infty}(S^n)$:
$$
C_{+,x_0}^{\infty}(S^n) = \{ u\in C_+^{\infty}(S^n) \ | \ u(x_0)= 1 \ \}.
$$ 
Then let ${\cal C}_{x_0}(S^n)$ be the orbit space of the induced
action of $C_{+,x_0}^{\infty}(S^n)$ on $\Riem(S^n)$. Clearly there is
a canonical projection map $p_1: {\cal C}_{x_0}(S^n) \lra {\cal
C}(S^n)$ which is a homotopy equivalence since $p_1^{-1}(C) \cong
\R$. Let $ {\cal C}_{x_0}^+(S^n) = p_1^{-1} \( {\cal C}^+(S^n)\).  $
We consider the following subgroup of the diffeomorphism group
$\Diff_+(S^n)$:
$$
\Diff_{x_0,+}(S^n)=\{\phi\in \Diff_+(S^n) \ | \ \phi(x_0)=x_0, \ \ d\phi=Id :
TM_{x_0} \rightarrow TM_{x_0} \}.
$$
The group $\Diff_{x_0,+}(S^n)$ acts (on the right, by pulling back a
metric) on the spaces ${\cal C}(S^n)$ and ${\cal C}_{x_0}(S^n) $.
Then it is an easy observation that the group $\Diff_{x_0,+}(S^n)$
acts freely on the space ${\cal C}_{x_0}(S^n)$.  Clearly the space
${\cal C}_{x_0}^+(S^n) $ of positive conformal classes is invariant
under this action. We define the moduli space of positive conformal
structures as the orbit space of the action of $\Diff_{x_0,+}(S^n)$ on
${\cal C}_{x_0}^+(S^n)$:
$$
{\cal M}_{x_0,\conf}^{+}(S^n) = {\cal C}_{x_0}^+(S^n)/\Diff_{x_0,+}(S^n).
$$
To make this construction usefull, we let
$\widetilde{\Diff}_{x_0,+}(S^n)\subset \Diff_+(S^n)$ be yet another
subgroup of diffeomorphisms $\phi$ with $\phi(x_0)=x_0$. The groups
$\Diff_{x_0,+}(S^n)$, $\widetilde{\Diff}_{x_0,+}(S^n)$, and
$\Diff_+(S^n)$ are clearly related to each other.\footnote{\ \ We are
grateful to Thomas Schick for a clarifying discussion on that
subject.} \ Indeed, one has the following fiber bundles:
$$
\begin{array}{l}
\Diff_{x_0,+}(S^n) \lra \widetilde{\Diff}_{x_0,+}(S^n) \lra GL^+(n; \R),
\\
\\
\widetilde{\Diff}_{x_0,+}(S^n) \lra \Diff_+(S^n) \lra S^n
\end{array}
$$
In particular, one concludes the isomorphisms:
$$
\pi_0 \Diff_{x_0,+}(S^n) \cong \pi_0 \widetilde{\Diff}_{x_0,+}(S^n)
\cong \pi_0 \Diff_+(S^n) \cong \Theta^{n+1},
$$ 
where $\Theta^{n+1}= \pi_0 \Diff_{+}(S^n)\cong \pi_0
\Diff_{x_0,+}(S^n)$ is the group of homotopy spheres.  The space
${\cal C}_{x_0}(S^n)$ is contractible, thus the orbit space ${\cal
C}_{x_0}(S^n)/ \Diff_{x_0,+}(S^n)$ is homotopy equivalent to the
classifying space $B\Diff_{x_0,+}(S^n)$. We obtain the following
commutative diagram of fiber bundles:
$$
\begin{diagram}
\setlength{\dgARROWLENGTH}{1.2em}
	\node{{\cal C}^+_{x_0}(S^n)} 
		\arrow[3]{e,t}{\subset}
		\arrow{s,l}{\Diff_{x_0,+}(S^n)}
	\node[3]{{\cal C}_{x_0}(S^n)} 
		\arrow{s,l}{\Diff_{x_0,+}(S^n)}
\\
	\node{{\cal M}_{x_0,\conf}^+(S^n)}
		\arrow[3]{e,t}{j}
	\node[3]{B\Diff_{x_0,+}(S^n)}
\end{diagram}
$$
In particular, one has the exact sequence in homotopy groups:
$$
\cdots\rightarrow \pi_1 ({\cal C}_{x_0}^+(S^n)) \stackrel{p_*}{\lra}
\pi_1 {\cal M}_{x_0,\conf}^+(S^n) \stackrel{\p}{\lra} \Theta^{n+1}
\stackrel{i_*}{\lra} \pi_0({\cal C}_{x_0}^+(S^n)) \stackrel{p_*}{\lra}
\pi_0 {\cal M}_{x_0,\conf}^+(S^n)
$$
We think that the moduli space ${\cal M}_{x_0,\conf}^+(S^n)$ is an
adequate model to study the positive scalar curvature metrics. It
captures all homotopy properties of the standard moduli space ${\cal
M}^{+}(S^n)$ of psc-metrics, and, on the other hand, is well-designed
for conformal geometry. We conclude with the following challenging problem.
\vspace{2mm}

\noindent
{\bf Problem 2.} {\sl Describe a rational homotopy type of the space ${\cal
M}_{x_0,\conf}^+(S^n)$.}

  
%
%
%

\end{document}